\documentclass[12pt, twoside]{article}
\usepackage{secdot}

\usepackage{amsmath, amsthm, amscd, amsfonts, amssymb, graphicx, color}
\usepackage[bookmarksnumbered, colorlinks, plainpages]{hyperref}

\begin{document}

\begin{center}
\Large{\bf $\ast$-Ricci-Yamabe Soliton and Contact Geometry}
\end{center}
\vspace{0.2cm}

\newcommand{\mvec}[1]{\mbox{\bfseries\itshape #1}}

\centerline{\large{Dibakar Dey}}
\vspace{0.2cm}
\centerline{Department of Pure Mathematics} 
\centerline{ University of Calcutta,
35, Ballygunge Circular Road}
\centerline{Kolkata-700019, West Bengal, India }
\centerline{E-mail: deydibakar3@gmail.com}



\newtheorem{theorem}{\quad Theorem}[section]

\newtheorem{definition}[theorem]{\quad Definition}

\newtheorem{proposition}[theorem]{\quad Proposition}

\newtheorem{question}[theorem]{\quad Question}

\newtheorem{remark}[theorem]{\quad Remark}

\newtheorem{corollary}[theorem]{\quad Corollary}

\newtheorem{note}[theorem]{\quad Note}

\newtheorem{lemma}[theorem]{\quad Lemma}
\newtheorem{example}[theorem]{\quad Example}
\newtheorem{notation}[theorem]{\quad Notation}
\numberwithin{equation}{section}
\newcommand{\be}{\begin{equation}}
\newcommand{\ee}{\end{equation}}
\newcommand{\bea}{\begin{eqnarray}}
\newcommand{\eea}{\end{eqnarray}}

\vspace{0.5 cm}
\textbf{Abstract:} It is well known that a unit sphere admits Sasakian 3-structure. Also, Sasakian manifolds are locally isometric to a unit sphere under several curvature and critical conditions. So, a natural question is: Does there exist any curvature or critical condition under which a Sasakian 3-manifold represents a geometrical object other than the unit sphere? In this regard,  as an extension of the $\ast$-Ricci soliton, the notion of $\ast$-Ricci-Yamabe soliton is introduced and studied on two classes contact metric manifolds. A $(2n + 1)$-dimensional non-Sasakian $N(k)$-contact metric manifold admitting $\ast$-Ricci-Yamabe soliton is completely classified. Further, it is proved that if a Sasakian 3-manifold $M$ admits $\ast$-Ricci-Yamabe soliton $(g,V,\lambda,\alpha,\beta)$ under certain conditions on the soliton vector field $V$, then $M$ is $\ast$-Ricci flat, positive Sasakian and the transverse geometry of $M$ is Fano. In addition, the Sasakian 3-metric $g$ is homothetic to a Berger sphere and the soliton is steady. Also, the potential vector field $V$ is an infinitesimal automorphism of the contact metric structure. \\



\section{Introduction }
\noindent In 2014, Kaimakamis and Panagiotidou \cite{kp} introduced
the notion of $\ast$-Ricci soliton from the Ricci soliton as \bea
\nonumber \frac{1}{2}\mathcal{L}_V g + S^\ast = \lambda g, \eea
where $\mathcal{L}_V$ denotes the Lie derivative along the vector field $V$, $S^\ast$ is the $\ast$-Ricci tensor given by $S^\ast(X,Y) = g(Q^\ast X,Y)$, $Q^\ast$ is the $\ast$-Ricci operator, $g$ is the Riemannian metric and $\lambda$ is a constant. The $\ast$-Ricci tensor $S^\ast$ is not symmetric in general. This definition is inconsistent if the $\ast$-Ricci tensor is not symmetric. For a consistent $\ast$-Ricci soliton on some space, the $\ast$-Ricci tensor of that space has to be symmetric. In \cite{gp}, Ghosh and Patra studied the notion of $\ast$-Ricci soliton on Sasakian manifolds. Further, in \cite{dd}, the authors studied $\ast$-Ricci soliton and $\ast$-gradient Ricci soliton within the framework of trans-Sasakian 3-manifolds.\\

\noindent In 2019, G$\ddot{u}$ler and Crasmareanu \cite{gc} proposed
the notion of Ricci-Yamabe flow on a Riemannian manifold $(M^n,g)$
as
 \bea
 \nonumber \frac{\partial g}{\partial t}(t) + 2\alpha S(t) + \beta r(t)g(t) = 0,
\eea
where $g$ is the Riemannian metric, $S$ is the Ricci tensor, $r$ is the scalar curvature and $\alpha$, $\beta$ are two real constants. Since $\alpha$ and $\beta$ are arbitrary constants, then one can freely choose the signs of $\alpha$ and $\beta$. This freedom of choice is very useful in differential geometry and theory of relativity. In \cite{akrami,bc}, the authors study space-time geometry with a bi-metric approach.\\

\noindent In \cite{ddey}, the present author proposed the idea of
Ricci-Yamabe soliton from the Ricci-Yamabe flow on a Riemannian
manifold $(M^n,g)$ as \bea \nonumber \mathcal{L}_V g + 2\alpha S =
(2\lambda - \beta r)g, \eea
where $\lambda,\; \alpha,\; \beta \in \mathbb{R}$. This generalizes a large class of soliton like equations such as Ricci soliton, Yamabe soliton, Einstein soliton, $\rho$-Einstein soliton. In \cite{ddey}, the present author classified two  classes of almost Kenmotsu manifolds admitting Ricci-Yamabe soliton.\\

\noindent We now extent the notion of $\ast$-Ricci soliton on a contact metric manifold to a more generalized version as follows:
\begin{definition}
A  contact metric manifold $(M,g)$ is said to admit a
$\ast$-Ricci-Yamabe soliton (in short, $\ast$-RYS)
$(g,V,\lambda,\alpha,\beta)$ if \bea \mathcal{L}_V g + 2\alpha
S^\ast = (2\lambda - \beta r^\ast)g, \label{1.1} \eea where
$\lambda,\;\alpha,\;\beta\; \in \; \mathbb{R}$ such that $\alpha
\neq 0$ and $r^\ast$ is the $\ast$-scalar curvature defined by
$r^\ast = trace(Q^\ast)$.
\end{definition}
\noindent If $V$ is gradient of some smooth function $f$ on $M$,
then it is called a $\ast$-gradient Ricci-Yamabe soliton (in short,
$\ast$-GRYS) and then \eqref{1.1} reduces to \bea \nabla^2 f +
\alpha S^\ast = (\lambda - \frac{1}{2}\beta r^\ast)g, \label{1.2}
\eea
where $\nabla^2 f$ is the Hessian of $f$ defined by $Hess_f(X,Y) = g(\nabla_X Df,Y)$, $D$ is the gradient operator. The $\ast$-RYS (or $\ast$-GRYS) is said to be expanding, steady or shrinking according as $\lambda$ is negative, zero or positive respectively. Note that, the $\ast$-RYS reduces to a $\ast$-Ricci soliton if $\alpha = 1$ and $\beta = 0$.\\

\noindent  $N(k)$-contact metric manifolds are special kind of
contact metric manifolds that generalizes Sasakian manifolds.
Sasakian geometry is an odd dimensional analog of the Kaehler
geometry and is an interesting topic to physicists as it perceived
relevance in string theory (see \cite{candelas,jm}). Due to
this connection with physics, we consider the notion of $\ast$-RYS
and $\ast$-GRYS within the frameworks of $N(k)$-contact geometry and
Sasakian 3-geometry. The paper is organized as follows: In section
2, we recall some basic relations and definitions on $N(k)$-contact
and Sasakian geometry. Section 3 deals with $N(k)$-contact metric as
a $\ast$-RYS. In section 4, a Sasakian 3-manifold is completely
classified admitting $\ast$-RYS under certain conditions on the soliton vector field.

\section{Preliminaries}
An odd dimensional smooth manifold $M$ together with a structure
$(\varphi, \xi, \eta,g)$ satisfying \be \varphi^{2}X = -X +
\eta(X)\xi,\; \eta(\xi)=1,\; \varphi \xi = 0,\; \eta \circ \varphi =
0, \label{2.1} \ee \be
 g(\varphi X,\varphi Y)=g(X,Y)-\eta(X)\eta(Y) \label{2.2}
\ee for all vector fields $X$, $Y$ on $M$ is called an almost
contact metric manifold, where $\varphi$ is a $(1,1)$-tensor field,
$\xi$ is a unit vector field called the characteristic vector field,
$\eta$ is a 1-form dual to $\xi$ and $g$ is the Riemannian metric.
It is easy to see from \eqref{2.2} that $\phi$ is skew-symmetric,
that is, \be g(\varphi X,Y) = - g(X,\varphi Y). \label{2.3} \ee
 A contact metric manifold is an almost contact metric manifold with $d\eta = g(X,\varphi Y)$. On a contact metric manifold, the $(1,1)$-tensor field $h$ is defined as $h = \frac{1}{2}\mathcal{L}_\xi \varphi$. The tensor field $h$ is symmetric and satisfies
 \be
 h\varphi = - \varphi h,\;\; trace(h) = trace(\varphi h) = 0,\;\; h\xi = 0. \label{2.4}
 \ee
 Also, on a contact metric manifold, we have
 \be
 \nabla_X \xi = - \varphi X - \varphi hX. \label{2.5}
 \ee
 In \cite{tanno}, Tanno introduced the notion of $k$-nullity distribution on a Riemannian manifold as
\bea \nonumber N(k) = \{ Z \in T(M) : R(X,Y)Z = k[g(Y,Z)X -
g(X,Z)Y]\}, \eea where $k \in \mathbb{R}$ and $T(M)$ is the Lie
algebra of all vector fields on $M$. If the characteristic vector
field $\xi \in N(k)$, then we call a contact metric manifold as
$N(k)$-contact metric manifold \cite{tanno}. For a $(2n +
1)$-dimensional $N(k)$-contact metric manifold, we have ( see
\cite{blair1,blair2}) \bea h^2 = (k - 1)\varphi^2,
\label{2.6} \eea \bea R(X,Y)\xi = k[\eta(Y)X - \eta(X)Y],
\label{2.7} \eea \bea R(\xi,X)Y = k[g(X,Y)\xi - \eta(Y)X],
\label{2.8} \eea \bea (\nabla_X \eta)Y = g(X + hX,\varphi Y),
\label{2.9} \eea \bea (\nabla_X \varphi)Y = g(X + hX,Y)\xi -
\eta(Y)(X + hX) \label{2.10} \eea
for all vector fields $X$, $Y$ on $M$, where $R$ is the Riemann curvatue tensor.\\

\noindent If the characteristic vector field $\xi$ is Killing type,
then a contact metric manifold is called a $K$-contact manifold and
if the structure $(\varphi,\xi,\eta,g)$ is normal, then a contact
metric manifold is called Sasakian. Also, an almost contact meteic
manifold is Sasakian if and only if \bea (\nabla_X \varphi Y) =
g(X,Y)\xi - \eta(Y)X \label{2.11} \eea holds for all vector fields
$X$, $Y$ on $M$. Hence, a $N(k)$-contact metric manifold reduces to
a Sasakian one if $h = 0$, that is, $k = 1$. On a Sasakian
3-manifold, the following relations are well known: \bea \nabla_X
\xi = -\varphi X, \label{2.12} \eea \bea (\nabla_X \eta)Y =
g(X,\varphi Y), \label{2.13} \eea \bea R(X,Y)\xi = \eta(Y)X -
\eta(X)Y, \label{2.14} \eea \bea R(\xi,X)Y = g(X,Y)\xi - \eta(Y)X.
\label{2.15} \eea

\noindent Since a 3-dimensional Riemannian manifold is conformally
flat, it's curvature tensor can be expressed as \bea \nonumber
R(X,Y)Z &=& [S(Y,Z)X - S(X,Z)Y + g(Y,Z)QX - g(X,Z)QY] \\ && -
\frac{r}{2}[g(Y,Z)X - g(X,Z)Y], \label{2.16} \eea where $r$ is the
scalar curvature and $S$ is the Ricci tensor. The Ricci tensor of a
Sasakian 3-manifold can be obtained from here as \bea S(X,Y) =
\frac{1}{2}[(r - 2)g(X,Y) + (6 - r)\eta(X)\eta(Y)]. \label{2.17}
\eea Note that the scalar curvature $r$ is not constant in general.
We now close this section by recalling the following definition:

\begin{definition}
Let $V$ be a vector field on an almost contact metric manifold $M$.
If there exist a smooth function $\sigma$ on $M$ such that
$\mathcal{L}_V g =  2\sigma g$, then $V$ is called a conformal
vector field. In particular, if $\sigma$ is constant then $V$ is homothetic and if $\sigma = 0$, then $V$ is called a
Killing vector field. Also $V$ is said to be an infinitesimal
contact transformation if $\mathcal{L}_V \eta = \psi \eta$ for some
smooth function $\psi$ on $M$. If $\psi = 0$, then $V$ is said to be
a strict infinitesimal contact transformation. If $V$ leaves all the
structure tensor fields $\varphi$, $\xi$, $\eta$ and $g$ invariant,
then $V$ is called an infinitesimal automorphism of the contact
metric structure.
\end{definition}

\section{$N(k)$-Contact Metric as a $\ast$-RYS}
In this section, we study the notion of $\ast$-RYS within the
framework of $N(k)$-contact metric manifolds. To prove the main
theorem of this section, we need the following lemmas:

\begin{lemma} (\cite{blair3}) \label{l3.1}
A contact metric manifold $M^{2n+1}$ satisfying the condition
$R(X,Y)\xi = 0$ for all $X$, $Y$ is locally isometric to the
Riemannian product of a flat $(n + 1)$-dimensional manifold and an
$n$-dimensional manifold of positive curvature $4$, i.e.,
$E^{n+1}(0) \times S^n(4)$ for $n > 1$ and flat for $n = 1$.
\end{lemma}

\begin{lemma} (\cite{dey2}) \label{l3.2}
A $(2n + 1)$-dimensional $N(k)$-contact metric manifold is
$\ast$-$\eta$-Einstein and the $\ast$-Ricci tensor is given by \bea
S^\ast (X,Y) = -k[g(X,Y) - \eta(X)\eta(Y)]. \label{3.1} \eea
\end{lemma}

\begin{note}
We observe that the $\ast$-Ricci tensor of a $N(k)$-contact metric
manifold is symmetric. Hence, the notion of $\ast$-RYS is consistent
in this setting. We are now ready to prove the main result of this
section.
\end{note}

\begin{theorem} \label{t3.4}
 If a $(2n + 1)$-dimensional non-Sasakian $N(k)$-contact metric manifold $M$ admits $\ast$-RYS $(g,V,\lambda,\alpha,\beta)$, then
 \begin{itemize}
 \item[$(1)$] The manifold $M$ is $\ast$-Ricci flat.
 \item[$(2)$] The manifold $M$ is locally isometric to $E^{n+1}(0) \times S^n(4)$ for $n > 1$ and flat for $n = 1$.
 \item[$(3)$] The soliton vector field V is homothetic.
 \end{itemize}
\end{theorem}

\begin{proof}
We start the proof by taking $g$-trace of \eqref{3.1} that gives
$r^\ast = -2nk$. Now, substituting \eqref{3.1} and the value of
$r^\ast$ in \eqref{1.1} yields
 \bea
 (\mathcal{L}_V g)(X,Y) = [2\lambda + 2nk\beta + 2k\alpha]g(X,Y) - 2k\alpha \eta(X)\eta(Y). \label{3.2}
 \eea
Differentiating the foregoing equation covariantly along any vector
field $Z$, we get \bea (\nabla_Z \mathcal{L}_V g)(X,Y) =
-2k\alpha[\eta(Y)(\nabla_Z \eta)X + \eta(X)(\nabla_Z \eta)Y].
\label{3.3} \eea Applying \eqref{2.9} in \eqref{3.3} yields \bea
(\nabla_Z \mathcal{L}_V g)(X,Y) = -2k\alpha[\eta(Y)g(Z + hZ,\varphi
X) + \eta(X)g(Z + hZ,\varphi Y)]. \label{3.4} \eea The well known
commutation formula (see \cite{yano}) \be \nonumber (\mathcal{L}_V
\nabla_X g - \nabla_X \mathcal{L}_V g - \nabla_{[V,X]}g)(Y,Z) =
-g((\mathcal{L}_V \nabla)(X,Y),Z) - g((\mathcal{L}_V
\nabla)(X,Z),Y). \ee leads to

\bea \nonumber g((\mathcal{L}_V \nabla)(X,Y),Z) &=&
\frac{1}{2}(\nabla_X \mathcal{L}_V g)(Y,Z) + \frac{1}{2}(\nabla_Y
\mathcal{L}_V g)(X,Z) \\ \nonumber && - \frac{1}{2}(\nabla_Z
\mathcal{L}_V g)(X,Y). \eea Using \eqref{3.4} in the preceding
equation yields \bea \nonumber g((\mathcal{L}_V \nabla)(X,Y),Z) &=&
2k\alpha[\eta(X)g(\varphi Y,Z)  + \eta(Y)g(\varphi X,Z) -
g(hX,\varphi Y)\eta(Z)], \eea which implies \bea
 (\mathcal{L}_V \nabla)(X,Y) = 2k\alpha[\eta(X)\varphi Y + \eta(Y)\varphi X - g(hX,\varphi Y)\xi]. \label{3.5}
\eea Substituting $Y = \xi$ in the foregoing equation, we get \bea
(\mathcal{L}_V \nabla)(X,\xi) = 2k\alpha \varphi X.\label{3.6} \eea
Now, \bea \nonumber (\nabla_Y \mathcal{L}_V \nabla)(X,\xi) =
\nabla_Y (\mathcal{L}_V \nabla)(X,\xi) - (\mathcal{L}_V
\nabla)(\nabla_Y X,\xi) - (\mathcal{L}_V \nabla)(X,\nabla_Y \xi).
\eea Using \eqref{3.5},\eqref{3.6},\eqref{2.5} and \eqref{2.10} in
the previous equation, we obtain \bea \nonumber (\nabla_Y
\mathcal{L}_V \nabla)(X,\xi) &=& 2k\alpha[g(Y + hY,X)\xi - \eta(X)(Y
- \eta(Y)\xi + hY)  \\ &&- \eta(X)(Y + hY) + g(hX,Y - \eta(Y)\xi +
hY)\xi] \label{3.7} \eea Due to Yano \cite{yano}, we have \bea
\nonumber (\mathcal{L}_V R)(X,Y)Z = (\nabla_X \mathcal{L}_V
\nabla)(Y,Z) - (\nabla_Y \mathcal{L}_V \nabla)(X,Z), \eea Using
\eqref{3.7} in the above formula, we obtain \bea \nonumber
(\mathcal{L}_V R)(X,\xi)\xi &=& (\nabla_X \mathcal{L}_V
\nabla)(\xi,\xi) - (\nabla_\xi \mathcal{L}_V \nabla)(X,\xi) \\ &=&
-4k\alpha[X - \eta(X)\xi + hX]. \label{3.8} \eea Now, replacing $Y$
by $\xi$ in \eqref{3.2} gives \bea \nonumber (\mathcal{L}_V
g)(X,\xi)  = [2\lambda + 2nk\beta]\eta(X), \eea which leads to \bea
(\mathcal{L}_V \eta)X - g(X,\mathcal{L}_V \xi) = [2\lambda +
2nk\beta]\eta(X). \label{3.9} \eea Putting $X = \xi$ in the
preceding equation, we can easily obtain that \bea
\eta(\mathcal{L}_V \xi) = -[\lambda + nk\beta]. \label{3.10} \eea
Now, with the help of \eqref{3.9}, \eqref{3.10},\eqref{2.7} and
\eqref{2.8}, we obtain \bea
 (\mathcal{L}_V R)(X,\xi)\xi = k[2\lambda + 2nk\beta](X - \eta(X)\xi). \label{3.11}
\eea Equating \eqref{3.8} and \eqref{3.11},  we have \bea k[2\lambda
+ 2nk\beta + 4\alpha](X - \eta(X)\xi) = -4k\alpha hX. \label{3.12}
\eea Taking $g$-trace of \eqref{3.12} and using \eqref{2.4} yields
\bea k[2\lambda + 2nk\beta + 4\alpha] = 0. \label{3.13} \eea
Operating $h$ on \eqref{3.12}  and using \eqref{2.6} leads to \bea
\nonumber k^2[(2\lambda + 2nk\beta + 4\alpha)^2 + 16\alpha^2(k -
1)](X - \eta(X)\xi) = 0, \eea which implies \bea k^2[(2\lambda +
2nk\beta + 4\alpha)^2 + 16\alpha^2(k - 1)] = 0. \label{3.14} \eea
Now, \eqref{3.13} and \eqref{3.14} together implies either $k = 0$
or the following relations holds: \bea 2\lambda + 2nk\beta + 4\alpha
= 0 \label{3.15} \eea and \bea (2\lambda + 2nk\beta + 4\alpha)^2 +
16\alpha^2(k - 1) = 0. \label{3.16} \eea
\textbf{Case 1:} If $k = 0$, then from (\ref{3.1}), we have $S^\ast = 0$, that is, the manifold is $\ast$-Ricci flat. Again from (\ref{2.7}), we have $R(X,Y)\xi = 0$ and hence, from lemma \ref{l3.1}, it follows that the manifold $M$ is locally isometric to $E^{n+1}(0) \times S^n(4)$ for $n > 1$ and flat for $n = 1$. Also, equation \eqref{3.2} reduces to $\mathcal{L}_V g = 2\lambda g$, where $\lambda$ is a constant. This shows that $V$ is homothetic.\\

\noindent \textbf{Case 2.} If $k \neq 0$, then \eqref{3.15} and
\eqref{3.16} holds. Since $\alpha \neq 0$, these two equations
implies that $k = 1$. This shows that $M$ reduces to a Sasakian
manifold, a contradiction to our hypothesis. This completes the
proof.
\end{proof}

\noindent A $\ast$-RYS reduces to a $\ast$-Ricci soliton if $\alpha
= 1$ and $\beta = 0$. Hence, from the above theorem, we can state
the following:
\begin{corollary}
If $(g,V,\lambda)$ is a $\ast$-Ricci soliton on a $(2n +
1)$-dimensional non-Sasakian $N(k)$-contact metric manifold $M$,
then
 \begin{itemize}
 \item[$(1)$] The manifold $M$ is $\ast$-Ricci flat.
 \item[$(2)$] The manifold $M$ is locally isometric to $E^{n+1}(0) \times S^n(4)$ for $n > 1$ and flat for $n = 1$.
 \item[$(3)$] The soliton vector field V is homothetic.
 \end{itemize}
\end{corollary}
\section{Sasakian 3-Metric as a $\ast$-RYS}
 In the preceding section, the Sasakian case is omitted. In this section, we consider an extended notion of $\ast$-RYS by considering the soliton vector field $V$ as a gradient vector field and an infinitesimal contact transformation on Sasakian 3-manifold. We start this section with the following discussion:\\

\noindent A contact metric manifold $M$ is said to be
$\eta$-Einstein if there exist two smooth functions $a$ and $b$ such
that
 \bea
 \nonumber S(X,Y) = ag(X,Y) + b\eta(X)\eta(Y)
 \eea
 for all vector fields $X$, $Y$ on $M$. A $(2n + 1)$-dimensional $\eta$-Einstein Sasakian manifold with $a = - 2$ and $b = 2n + 2$ is known as null-Sasakian. An example of such a manifold is a Sasakian space form $\mathbb{R}^{2n+1}$ with constant $\varphi$-sectional curvature $-3$, which is identifiable with a Heisenberg group. Also, $\eta$-Einstein Sasakian manifold with $a > - 2$ is called positive-Sasakian. In this case, the transverse geometry of $M$ is Fano. For more details, we refer the reader to go through \cite{cp2}. From \eqref{2.17}, we see that a Sasakian 3-manifold is $\eta$-Einstein. Thus a Sasakian 3-manifold is null-Sasakian if $r = - 2$ and  positive-Sasakian if $r > - 2$. For $r > - 2$, the transverse geometry of the Sasakian 3-manifold is Fano.\\

 \noindent In \cite{gp}, Ghosh and Patra obtained the expression of the $\ast$-Ricci tensor for a $(2n + 1)$-dimensional Sasakian manifold which involves the Ricci tensor. In \cite{dd}, the present author together with Majhi presents a expresion of the $\ast$-Ricci tensor for a trans-Sasakian 3-manifold, which is a generalization of the Sasakian 3-manifold. From lemma 3.1 of \cite{dd} or using \eqref{2.17} in lemma 3.1 of \cite{gp}, we obtain the expression of the $\ast$-Ricci tensor for Sasakian 3-manifold $M$ as
\bea S^\ast (X,Y) = \frac{1}{2}(r - 4)[g(X,Y) - \eta(X)\eta(Y)]
\label{4.1} \eea for all vector fields $X$, $Y$ on $M$. Note that,
$S^\ast$ is symmetric. Taking $g$-trace of \eqref{4.1}, we have
$r^\ast = r - 4$. We now consider $V$ as gradient of some smooth
function $f : M \rightarrow \mathbb{R}$ in the definition of
$\ast$-RYS and prove the following theorem.

\begin{theorem} \label{t4.1}
If $(g,V,\lambda,\alpha,\beta)$ be a $\ast$-RYS on a Sasakian
3-manifold $M$ such that either $(i)$ $V$ is a gradient vector field or $(ii)$ $V$ is an infinitesimal contact transformation, then
\begin{itemize}
     \item[$(1)$] The manifold $M$ is $\ast$-Rici flat.
     \item[$(2)$] The Sasakian 3-manifold $M$ is positive-Sasakian and the transverse geometry of $M$ is Fano.
     \item[$(3)$] The Sasakian 3-metric $g$ is homothetic to a Berger sphere.
     \item[$(4)$] The soliton vector field $V$ is an infinitesimal automorphism of the contact metric structure.
     \item[$(5)$] The $\ast$-RYS is steady.
 \end{itemize}
\end{theorem}

\begin{proof}
$(i)$ If the soliton vector field $V$ is gradient of some smooth function $f$ on $M$, then the equation \eqref{1.2} can be exhibited as \bea \nabla_X Df =
[\lambda - \frac{1}{2}(r - 4)\beta]X - \alpha Q^\ast X \label{4.2}
\eea for any vector field $X$ on $M$, where $V = Df$ and $r^\ast = r
- 4$ is used. Differentiating \eqref{4.2} covariantly along any
vector field $Y$ yields \bea
 \nabla_Y \nabla_X Df = [\lambda - \frac{1}{2}(r - 4)\beta]\nabla_Y X  - \frac{1}{2}\beta(Yr)X -  \alpha \nabla_Y Q^\ast X. \label{4.3}
\eea Interchanging $X$ and $Y$ in the preceding equation, we have
\bea
 \nabla_X \nabla_Y Df = [\lambda - \frac{1}{2}(r - 4)\beta]\nabla_X Y  - \frac{1}{2}\beta(Xr)Y -  \alpha \nabla_X Q^\ast Y. \label{4.4}
\eea From \eqref{4.2}, we get \bea \nabla_{[X,Y]} Df = [\lambda -
\frac{1}{2}(r - 4)\beta](\nabla_X Y - \nabla_Y X) - \alpha Q^\ast
(\nabla_X Y - \nabla_Y X). \label{4.5} \eea Now, the curvature
tensor $R$ is given by 
\bea \nonumber R(X,Y)Df = \nabla_X \nabla_Y Df - \nabla_Y \nabla_X Df
- \nabla_{[X,Y]} Df. 
\eea 
Using \eqref{4.3}-\eqref{4.5} in the
foregoing equation, we get \bea
 R(X,Y)Df =  \frac{1}{2}\beta(Yr)X - \frac{1}{2}\beta(Xr)Y - \alpha[(\nabla_X Q^\ast)Y - (\nabla_Y Q^\ast)X]. \label{4.6}
\eea From \eqref{4.1}, we have \bea Q^\ast X = \frac{1}{2}(r - 4)[X
- \eta(X)\xi]. \label{4.7} \eea With the help of \eqref{2.12} and
\eqref{2.13}, from \eqref{4.7}, we obtain \bea
 (\nabla_Y Q^\ast)X = \frac{1}{2}(Yr)[X - \eta(X)\xi] - \frac{1}{2}(r - 4)[g(\varphi X,Y)\xi - \eta(X)\varphi Y]. \label{4.8}
\eea Applying \eqref{4.8} in \eqref{4.6}, we obtain \bea \nonumber
R(X,Y)Df &=& \frac{1}{2}\beta[(Yr)X - (Xr)Y] - \frac{1}{2}\alpha
(Xr)[Y - \eta(Y)\xi] \\ \nonumber && + \frac{1}{2}\alpha (Yr)[X -
\eta(X)\xi] +  \frac{1}{2}\alpha (r - 4)[2g(X,\varphi Y)\xi \\ && -
\eta(Y)\varphi X + \eta(X)\varphi Y]. \label{4.9} \eea Substituting
$X = \xi$ in \eqref{4.9} and noting $(\xi r) = 0$ (as $\xi$ is
Killing), we infer that \bea \nonumber R(\xi,Y)Df =
\frac{1}{2}\beta(Yr)\xi + \frac{1}{2}\alpha (r - 4)\varphi Y. \eea
Taking inner product of the preceding equation with $X$ yields
 \bea
g(R(\xi,Y)Df,X) = \frac{1}{2}\beta(Yr)\eta(X) + \frac{1}{2}\alpha (r
- 4)g(\varphi Y,X). \label{4.10} \eea Since $g(R(\xi,Y)Df,X) = -
g(R(\xi,Y)X,Df)$, then using \eqref{2.15}, we get \bea
g(R(\xi,Y)Df,X) = - g(X,Y)(\xi f) + \eta(X)(Yf). \label{4.11} \eea
Equating \eqref{4.10} and \eqref{4.11} and then antisymmetrizing
yields \bea \nonumber && \frac{1}{2}\beta(Yr)\eta(X) -
\frac{1}{2}\beta(Xr)\eta(Y) + \alpha (r - 4)g(\varphi Y,X) \\ && =
\eta(X)(Yf) - \eta(Y)(Xf).  \label{4.12} \eea Replacing $X$ by $\phi
X$ and $Y$ by $\phi Y$ in \eqref{4.12} yields \bea \nonumber
\alpha(r - 4)g(\phi X,Y) = 0, \eea
which implies $r = 4$ as $\alpha \neq 0$ by definition of $\ast$-RYS. Hence, from \eqref{4.1}, we get $S^\ast = 0$. This proves (1).\\

\noindent Since $r = 4 > - 2$, then $M$ is positive-Sasakian and the transverse geometry of $M$ is Fano proving (2).\\

\noindent The Tanaka-Webster curvature (see \cite{bsg1}) of a Sasakian 3-manifold is given by $ W = \frac{1}{4}(r + 2)$. Since $r = 4$, then $W = \frac{3}{2}$. Following the classification given by Guilfoyle \cite{bsg2} for $0 < W < 2$, we conclude that $g$ is homothetic to a Berger sphere. This proves (3).\\

\noindent Now, equation \eqref{4.12} reduces to \bea \nonumber
\eta(X)(Yf) - \eta(Y)(Xf) = 0. \eea Putting $X = \xi$ in the
preceding equation yields \bea \nonumber (Yf) - \eta(Y)(\xi f) = 0,
\eea which implies $Df = (\xi f)\xi$. Therefore, we get $V = Df =
(\xi f)\xi$, that is, $V$ is pointwise collinear with $\xi$. For
simplicity, we write $c = (\xi f)$. Then \bea (\mathcal{L}_V g)(X,Y)
= (\mathcal{L}_{c \xi} g)(X,Y) = (Xc)\eta(Y) + (Yc)\eta(X).
\label{4.13} \eea Using \eqref{4.13}, $S^\ast = 0$ and $r^\ast = 0$
in \eqref{1.1}, we have \bea (Xc)\eta(Y) + (Yc)\eta(X) = 2\lambda
g(X,Y). \label{4.14} \eea Substituting $X = Y = \xi$ in \eqref{4.14}
yields \bea 2(\xi c) = 2\lambda. \label{4.15} \eea Let $\{e_i\}$ be
any orthonormal frame in $M$. Now, substituting $X = Y = e_i$ in
\eqref{4.14} and summing over $i$, we obtain \bea 2(\xi c) =
6\lambda. \label{4.16} \eea
The equtions \eqref{4.15} and \eqref{4.16}  together implies $\lambda = 0$ and hence, the $\ast$-GRYS is steady proving (5).\\

\noindent Thus we get $(\xi c) = 0$. Then using $\lambda = 0$ and
putting $Y = \xi$ in \eqref{4.14}, we get $(Xc) = 0$ for any vector
field $X$, which implies $c = (\xi f)$ is a constant. Therefore, $V$
is a constant multiple of $\xi$. Since $c$
is a constant, then \eqref{4.13} gives $\mathcal{L}_V g = 0$. It is
easy to see that $\mathcal{L}_V \xi = \mathcal{L}_{c\xi} \xi = 0$.
Now, $\mathcal{L}_V g = 0$ and $\mathcal{L}_V \xi = 0$ together
implies $\mathcal{L}_V \eta = 0$. With the help of \eqref{2.11}  and
\eqref{2.12}, it can be easily proved that $\mathcal{L}_V \varphi =
0$. Therefore, $V$ leaves all the structure tensors $\varphi$,
$\xi$, $\eta$ and $g$ invariant, that is, $V$ is an infinitesimal
automorphism of the contact metric structure proving (4) and this completes the proof of $(i)$.\\

$(ii)$ If $V$ is an infinitesimal contact transformation, then there exist a smooth function $f$ on $M$ such that
\bea
\mathcal{L}_V \eta = f\eta. \label{4.17}
\eea
Since $d\eta(X,Y) = g(X,\varphi Y)$, then
 \bea
\nonumber (\mathcal{L}_V d\eta)(X,Y) &=& \mathcal{L}_V d\eta(X,Y) - d\eta(\mathcal{L}_V X,Y) - d\eta(X,\mathcal{L}_V Y) \\ \nonumber &=& \mathcal{L}_V g(X,\varphi Y) - g(\mathcal{L}_V X,\varphi Y) - g(X,\varphi \mathcal{L}_V Y) \\ &=& (\mathcal{L}_V g)(X,\varphi Y) + g(X, (\mathcal{L}_V \varphi)Y). \label{4.18}
 \eea
Substituting the value of $(\mathcal{L}_V g)(X,\varphi Y)$ from \eqref{1.1} in \eqref{4.18}, we get
\be
(\mathcal{L}_V d\eta)(X,Y) = - 2\alpha S^\ast(X,\varphi Y) +  (2\lambda - \beta r^\ast)g(X,\varphi Y) + g(X,(\mathcal{L}_V \varphi)Y). \label{4.19}
\ee
Since \eqref{4.17} holds, then we have
\bea
 \mathcal{L}_V d\eta = d\mathcal{L}_V \eta = df \wedge \eta + fd\eta, \label{4.20}
\eea
 which implies
 \bea
 (\mathcal{L}_V d\eta)(X,Y) = \frac{1}{2}[(Xf)\eta(Y) - (Yf)\eta(X)] + fg(X,\varphi Y). \label{4.21}
 \eea
Equating \eqref{4.19} and \eqref{4.21}, we get
\bea
\nonumber g(X, (\mathcal{L}_V \varphi)Y) &=& \frac{1}{2}[(Xf)\eta(Y) - (Yf)\eta(X)] \\ \nonumber && + [f + \alpha(r - 4) - 2\lambda + \beta(r - 4)]g(X,\varphi Y),
\eea
which implies
\bea
(\mathcal{L}_V \varphi)Y = \frac{1}{2}[\eta(Y)Df - (Yf)\xi] + [f + \alpha(r - 4) - 2\lambda + \beta(r - 4)]\varphi Y. \label{4.22}
\eea
Replacing $Y$ by $\xi$ in the preceding equation, we have
 \bea
 (\mathcal{L}_V \varphi)\xi = \frac{1}{2}[Df - (\xi f)\xi]. \label{4.23}
 \eea
 Taking $g$-trace of the equation \eqref{1.1} yields
 \bea
 \operatorname{Div} V = \frac{3}{2}[2\lambda - \beta(r - 4)] - \alpha(r - 4), \label{4.24}
 \eea
 where `$\operatorname{Div}$' stands for divergence. Let $\Sigma$ be the volume form of $M$, that is, $\Sigma = \eta \wedge (d\eta)^n \neq 0$. Lie differentiating this along the vector field $V$ and applying the formula $\mathcal{L}_V \Sigma = (\operatorname{Div}V)\Sigma$ , \eqref{4.17} and \eqref{4.20}, we obtain $(\operatorname{Div}V)\Sigma = (n + 1)f\Sigma$, which implies
 \bea
 \operatorname{Div}V = (n + 1)f. \label{4.25}
 \eea
Equating \eqref{4.24} and \eqref{4.25}, we obtain
\bea
(n + 1)f = \frac{3}{2}[2\lambda - \beta(r - 4)] - \alpha(r - 4). \label{4.26}
\eea
From \eqref{1.1}, we write
\bea
\nonumber (\mathcal{L}_V g)(X,\xi) = [2\lambda - \beta(r - 4)]\eta(X),
\eea
 which implies
 \bea
 (\mathcal{L}_V \eta)X - g(X,\mathcal{L}_V \xi) = [2\lambda - \beta(r - 4)]\eta(X). \label{4.27}
 \eea
Using \eqref{4.17} in \eqref{4.27} gives
\bea
 f\eta(X) - g(X,\mathcal{L}_V \xi) = [2\lambda - \beta(r - 4)]\eta(X). \label{4.28}
 \eea
Substituting $X = \xi$ in the foregoing equation gives
\bea
\eta(\mathcal{L}_V \xi) = [f - 2\lambda + \beta(r - 4)]. \label{4.29}
\eea 
Now, substitution of $X = \xi$ in \eqref{4.27} yields
\bea
\eta(\mathcal{L}_V \xi) = -\frac{1}{2}[2\lambda - \beta(r - 4)]. \label{4.30}
\eea
Equating \eqref{4.29} and \eqref{4.30}, we get
\bea
f = \frac{1}{2}[2\lambda - \beta(r - 4)]. \label{4.31}
\eea
We now use the equation \eqref{4.31} in \eqref{4.28} to obtain
\bea
\mathcal{L}_V \xi = - f\xi. \label{4.32}
\eea
With the help of \eqref{4.32}, it can be easily seen that $(\mathcal{L}_V \varphi)\xi = 0$ and hence, equation \eqref{4.23} reduces to $Df = (\xi f)\xi$. Taking inner product of this with $Y$, we have
\bea
df(Y) = (\xi f)\eta(Y), \label{4.33}
\eea
which implies $df = (\xi f)\eta$. Now, taking exterior derivative of this, we have $d^2 f = d(\xi f) \wedge \eta + (\xi f)d\eta$. Applying $d^2 = 0$ on it and then taking wedge product with $\eta$ yields $(\xi f)\eta \wedge d\eta = 0$. Since $\eta \wedge d\eta \neq 0$ on $M$, then $(\xi f) = 0$ and therefore, equation \eqref{4.33} implies $f$ is constant. Integrating both sides of
\eqref{4.25} over $M$ and applying the divergence theorem we get $f = 0$. Thus, equations \eqref{4.17} and \eqref{4.32} implies $\mathcal{L}_V \eta = \mathcal{L}_V \xi = 0$. From \eqref{4.31}, we have $[2\lambda - \beta(r - 4)] = 0$ and hence, from \eqref{4.26}, we get $r = 4$ as $\alpha$ is a non-zero constant, which implies $\lambda = 0$ proving (5). Therefore, equation \eqref{4.22} gives $\mathcal{L}_V \varphi = 0$ and equation \eqref{4.1} gives $S^\ast = 0$. Now, from \eqref{1.1} yields $\mathcal{L}_V g = 0$. Thus $V$ is an infinitesimal automorphism of the contact metric structure. Since $r = 4$, then by similar arguments as in part $(i)$, the statements (1), (2) and (3) holds. This completes the theorem.
\end{proof}

 \begin{remark}
 If $M$ is complete and since $r = 4 > 0$, then by Myers theorem \cite{myers}, $M$ is necessarily compact.
 \end{remark}

\begin{remark}
The $\ast$-RYS reduces to a $\ast$-Ricci soliton if
$\alpha = 1$ and $\beta = 0$. Thus the results of the above theorem
holds  for a Sasakian 3-manifold admitting $\ast$-Ricci
soliton.
\end{remark}

\noindent It is proved that a $N(k)$-contact metric or a Sasakian 3-metric satisfying the $\ast$-RYS is $\ast$-Ricci flat. Thus a natural question is
\begin{question}
Does there exist an almost contact metric $g$ such that $g$ admits
$\ast$-RYS or $\ast$-GRYS whose $\ast$-Ricci tensor is not
identically zero or in other words, the manifold is not $\ast$-Ricci
flat ?
\end{question}

\subsection*{Further Scope of Study}
The notion of $\ast$-RYS is a generalized version of $\ast$-Ricci soliton. For instance, one can see that the equation \eqref{1.1} provides several soliton like equations such as
\begin{itemize}
\item $\ast$-Yamabe soliton for $\alpha = 0$ and $\beta = 2$.
\item $\ast$-Einstein soliton for $\alpha = 1$ and $\beta = -1$.
\item $\ast$-$\rho$-Einstein soliton for $\alpha = 1$ and $\beta = -2\rho$.
\end{itemize}
The notion of  $\ast$-Yamabe soliton is not introduced yet. So, this can be a platform for introducing the notion of $\ast$-Yamabe soliton. In this article, the author studied the notion of $\ast$-RYS within the framework of $N(k)$-contact and Sasakian geometry. There is a large class of almost contact metric manifolds to study this notion. Further, one can consider this new notion on para-contact manifolds or Lorentzian manifolds.

\end{document}